\documentclass{amsart}
\usepackage{amssymb} 
\usepackage{color}
\usepackage[all]{xy}
\usepackage{enumerate}
\usepackage{ulem}

 \newtheorem{theorem}{Theorem}[section]

 \newtheorem{proposition}[theorem]{Proposition}
 \newtheorem{lemma}[theorem]{Lemma}
 \newtheorem{corollary}[theorem]{Corollary}
 \newtheorem{remark}[theorem]{Remark}
 \newtheorem{definition}[theorem]{Definition}
 \newtheorem{notation}[theorem]{Notation}

 \numberwithin{equation}{section}
 \def\subrel#1#2{\mathrel{\mathop{#2}\limits_{#1}}}

\def\cat{{\rm Cat}_{\infty}}

\def\wcat{\widehat{\rm Cat}_{\infty}}

\begin{document}

\title
{Uniqueness of monoidal adjunctions}
\author{Takeshi Torii}
\address{Department of Mathematics, 
Okayama University,
Okayama 700--8530, Japan}
\email{torii@math.okayama-u.ac.jp}


\subjclass[2020]{18N70 (primary), 18N60, 55U40 (secondary)}
\keywords{Monoidal $\infty$-category, lax monoidal functor,
Day convolution, 
$\infty$-operad.}

\date{March 11, 2023 (version~2.0)}

\begin{abstract}
There are two dual equivalences
between the $\infty$-category of
$\mathcal{O}$-monoidal $\infty$-categories with
right adjoint lax $\mathcal{O}$-monoidal functors and
that with left adjoint oplax $\mathcal{O}$-monoidal functors,
where $\mathcal{O}$ is an $\infty$-operad.
We study the space of equivalences
between these two $\infty$-categories,
and show that the two equivalences
equipped with compatible $\mathcal{O}$-monoidal presheaf functors
are canonically equivalent.
\end{abstract}

\maketitle

\section{Introduction}

In this paper we study the space of 
equivalences between the $\infty$-category
of $\mathcal{O}$-monoidal $\infty$-categories
with right adjoint lax $\mathcal{O}$-monoidal functors
and that with left adjoint oplax $\mathcal{O}$-monoidal functors,
where $\mathcal{O}$ is an $\infty$-operad.
The purpose of this paper is to show a uniqueness result
of such equivalences 
equipped with compatible $\mathcal{O}$-monoidal presheaf functors.

In higher category theory
we expect that statements in classical category theory
should be generalized appropriately.
We say that an $\infty$-category in which $j$-morphisms are all
invertible for $j>r$
is an $(\infty,r)$-category.
The fact that $(\infty,1)$-categories
have non-invertible $1$-morphisms is essential
and it implies that $(\infty,1)$-categories
have similar categorical structure to
ordinary $1$-categories.
Therefore,
we expect that statements in classical $1$-category
should also hold in $(\infty,1)$-categories
by formulating them appropriately.

There are a lot of models for $(\infty,1)$-categories.
Complete Segal spaces,
Segal categories,
topological categories,
simplicial categories and quasi-categories
are models for $(\infty,1)$-categories.
It is known that these models are all equivalent 
(see, for example, \cite{Bergner,Joyal-Tierney}).
Quasi-categories are easy to handle
since they are closed under categorical constructions.
Quasi-categories were introduced
by Boardman-Vogt~\cite{Boardman-Vogt}
by name of weak Kan complexes
in their work on homotopy invariant algebraic structures.
The theory of quasi-categories
is further developed by Joyal~\cite{Joyal1}
and Lurie~\cite{Lurie1,Lurie2}.

In classical category theory
the left adjoint of a lax monoidal functor
between monoidal categories
is canonically an oplax monoidal functor
by doctrinal adjunction~\cite{Kelly}.
Dually,
the right adjoint of an oplax monoidal functor
is lax monoidal.
These correspondences give rise to a dual equivalence
between the category of monoidal categories
with right adjoint lax monoidal functors
and that with left adjoint oplax monoidal functors.

Generalizations of these correspondences
have been considered
in the setting of higher categories.
Lurie~\cite{Lurie2} showed that the right adjoint 
of a strong monoidal functor between monoidal $\infty$-categories
is lax monoidal.
This result is useful in many situations
but we would like to have more general correspondences
for monoidal adjunctions.
In \cite{Torii2}
we showed that the right adjoint 
of an oplax monoidal functor is lax monoidal,
and vice versa.
Haugseng, Hebestreit, Linskens, Nuiten~\cite{HHLN1}
formulated the correspondences in more precise way.
They constructed a bidual equivalence 
between the $(\infty,2)$-category
of $\mathcal{O}$-monoidal $\infty$-categories
with left adjoint oplax $\mathcal{O}$-monoidal
functors and that with right adjoint
lax $\mathcal{O}$-monoidal functors
for any $\infty$-operad $\mathcal{O}$
by introducing notions of curved orthofibrations
and (op-)Gray fibrations.
In \cite{Torii4},
on the other hand,
we proved a similar equivalence
restricted to the underlying $(\infty,1)$-categories
by another method.
We showed that there is a dual equivalence 
between 
the $\infty$-category 
${\rm Mon}_{\mathcal{O}}^{\rm lax, R}(\cat)$
of $\mathcal{O}$-monoidal
$\infty$-categories with right adjoint
lax $\mathcal{O}$-monoidal functors
and the $\infty$-category
${\rm Mon}_{\mathcal{O}}^{\rm oplax,L}(\cat)$
of $\mathcal{O}$-monoidal
$\infty$-categories with left adjoint
oplax $\mathcal{O}$-monoidal functors
by constructing a perfect pairing between them.

Thus,
we have two dual equivalences
between
${\rm Mon}_{\mathcal{O}}^{\rm lax, R}(\cat)$
and
${\rm Mon}_{\mathcal{O}}^{\rm oplax, L}(\cat)$.
One is the restriction of the equivalence 
obtained by 
Haugseng, Hebestreit, Linskens, Nuiten~\cite{HHLN1},
and the other is by
the author~\cite{Torii4}.
The purpose of this paper
is to give a uniqueness result of such equivalences.

In \cite{Torii4}
we have studied $\mathcal{O}$-monoidal presheaf functors
$\mathbb{P}_{\mathcal{O}!}:
           {\rm Mon}_{\mathcal{O}}^{\rm lax}(\cat)
\to
           {\rm Mon}_{\mathcal{O}}^{\rm lax}({\rm Pr}^{\rm L})$
           and
$\mathbb{P}_{\mathcal{O}}^*:
           {\rm Mon}_{\mathcal{O}}^{\rm oplax}(\cat)^{\rm op}
\to
           {\rm Mon}_{\mathcal{O}}^{\rm lax}({\rm Pr}^{\rm L})$.
In this paper
we consider a space ${\mathcal Eq}$
of equivalences $F$ between 
${\rm Mon}_{\mathcal{O}}^{\rm lax, R}(\cat)$
and
${\rm Mon}_{\mathcal{O}}^{\rm oplax, L}(\cat)^{\rm op}$
equipped with an equivalence
$\mathbb{P}_{\mathcal{O!}}\stackrel{\simeq}{\to}
\mathbb{P}_{\mathcal{O}}^*\circ F$.
The first main result of this paper
is the following theorem,
which means a uniqueness of such equivalences.

\begin{theorem}[{cf.~Theorem~\ref{thm:uniqueness-equivalences}}]
\label{theorem:main-theorem-1}
The space ${\mathcal Eq}$ is contractible.
\end{theorem}

By construction, 
the equivalence $T$ in \cite{Torii4}
is equipped with an equivalence
$\theta: \mathbb{P}_{\mathcal{O}!}\stackrel{\simeq}{\to}
\mathbb{P}_{\mathcal{O}}^*\circ T$.
Hence, we can regard the pair $(T,\theta)$
as an object of ${\mathcal Eq}$. 
For the restriction $H$
of the equivalence in \cite{HHLN1},
we will show that there is 
a canonical equivalence
$\gamma: \mathbb{P}_{\mathcal{O}!}\stackrel{\simeq}{\to}
\mathbb{P}_{\mathcal{O}}^*\circ H$
which is compatible with Yoneda embeddings.
Hence,
we can also regard the pair $(H,\gamma)$
as an object of ${\mathcal Eq}$.
By Theorem~\ref{theorem:main-theorem-1},
we obtain the following theorem,
which is the second main theorem of this paper.

\begin{theorem}[{cf.~Theorem~\ref{thm:T-eq-H}}]
\label{theorem:main-theorem-2}
The pair $(T,\theta)$ 
is canonically equivalent to $(H,\gamma)$.
\end{theorem}

This theorem implies the following corollary.

\begin{corollary}[{cf. Corollary~\ref{cor:main-claim}}]
The equivalence $T$ in \cite{Torii4}
is equivalent to the restriction $H$
of the equivalence in \cite{HHLN1}.
\end{corollary}

The organization of this paper is as follows:
In \S\ref{section:monoidal-presheaf-functors}
we review properties of monoidal presheaf functors
$\mathbb{P}_{\mathcal{O}!}$ and $\mathbb{P}_{\mathcal{O}}^*$
studied in \cite{Torii4}.
In \S\ref{section:uniqueness}
we prove a uniqueness result
(Theorem~\ref{theorem:main-theorem-1}).
We introduce the space ${\mathcal Eq}$
of equivalences between
${\rm Mon}_{\mathcal{O}}^{\rm lax,R}(\cat)$
and
${\rm Mon}_{\mathcal{O}}^{\rm oplax,L}(\cat)^{\rm op}$
compatible with monoidal presheaf functors,
and show that ${\mathcal Eq}$ is contractible.
In \S\ref{section:canonical-equivalence}
we prove the second main theorem
(Theorem~\ref{theorem:main-theorem-2}).
Showing that the restriction $H$
of the equivalence in \cite{HHLN1}
extends to an object $(H,\gamma)$
in ${\mathcal Eq}$,
we prove that the pair $(T,\theta)$
is canonically equivalent to $(H,\gamma)$.

\begin{notation}\rm
We fix Grothendieck universes 
$\mathcal{U}\in \mathcal{V}\in\mathcal{W}
\in \mathcal{X}$
throughout this paper. 
We say that an element of $\mathcal{U}$ is small,
an element of $\mathcal{V}$ is large, 
an element of $\mathcal{W}$ is very large,
and an element of $\mathcal{X}$
is super large.

We denote by ${\rm Cat}_{\infty}$
the large $\infty$-category of small $\infty$-categories,
and by $\widehat{\rm Cat}_{\infty}$
the very large $\infty$-category of large $\infty$-categories.
We denote by $\mathrm{Pr}^{\rm L}$ 
the subcategory of $\widehat{\rm Cat}_{\infty}$
spanned by presentable $\infty$-categories
and colimit-preserving functors.
We write $\mathcal{S}$
for the large $\infty$-category of small spaces,
$\widehat{\mathcal{S}}$
for the very large $\infty$-category
of large spaces,
and $\widetilde{S}$
for the super large $\infty$-category
of very large spaces.

For an $\infty$-category $\mathcal{C}$,
we denote by $h\mathcal{C}$
its homotopy category,
and by $\pi_{\mathcal{C}}: \mathcal{C}\to h\mathcal{C}$
the canonical functor.
For a functor $F: \mathcal{C}\to\mathcal{D}$
between $\infty$-categories,
we write $hF: h\mathcal{C}\to h\mathcal{D}$
for the induced functor on homotopy categories.

For an $\infty$-operad $\mathcal{O}^{\otimes}\to
{\rm Fin}_*$,
we write $\mathcal{O}$ for $\mathcal{O}^{\otimes}_{\langle 1\rangle}$.
For an $\mathcal{O}$-monoidal $\infty$-category
$\mathcal{C}^{\otimes}\to \mathcal{O}^{\otimes}$,
we denote by $\mathcal{C}$ 
the underlying $\infty$-category
$\mathcal{C}^{\otimes}\times_{\mathcal{O}^{\otimes}}\mathcal{O}$,
and we say that $\mathcal{C}$
is an $\mathcal{O}$-monoidal $\infty$-category
for simplicity.
\end{notation}

\noindent
{\bf Acknowledgements}.
The author would like to thank
Jonathan Beardsley
for letting him know 
references for the result by
Haugseng-Hebestreit-Linskens-Nuiten.
He would also like to thank the referee 
of the preprint~\cite{Torii4}
because this paper grows out of his trying
to answer questions and comments
by the referee.
The author was partially supported 
by JSPS KAKENHI Grant Numbers JP17K05253.

\section{Monoidal presheaf functors}
\label{section:monoidal-presheaf-functors}

In \cite{Torii4}
we studied monoidal
presheaf functors $\mathbb{P}_{\mathcal{O}!}$
and $\mathbb{P}_{\mathcal{O}}^*$.
In this section we review
the properties of
$\mathbb{P}_{\mathcal{O}!}$ and $\mathbb{P}_{\mathcal{O}}^*$
for reader's convenience.

For a small $\infty$-category $\mathcal{C}$,
we denote by ${\rm P}(\mathcal{C})$
the $\infty$-category 
${\rm Fun}(\mathcal{C}^{\rm op},\mathcal{S})$
of presheaves on $\mathcal{C}$
with values in $\mathcal{S}$,
and by $J_{\mathcal{C}}:
\mathcal{C}\to {\rm P}(\mathcal{C})$
the Yoneda embedding.
We have presheaf functors
${\rm P}_!: \cat\to {\rm Pr}^{\rm L}$ and
${\rm P}^*: \cat^{\rm op}\to {\rm Pr}^{\rm L}$,
where ${\rm P}_!(\mathcal{C})=
{\rm P}(\mathcal{C})={\rm P}^*(\mathcal{C})$.
For a functor $f:\mathcal{C}\to\mathcal{D}$,
we have 
${\rm P}_!(f)=f_!: {\rm P}(\mathcal{C})\to
{\rm P}(\mathcal{D})$
which is a left Kan extension of
$J_{\mathcal{D}}\circ f$ along $J_{\mathcal{C}}$,
and ${\rm P}^*(f)=f^*: {\rm P}(\mathcal{D})\to
{\rm P}(\mathcal{C})$
which is given by
the composition with $f^{\rm op}$.

In \cite{Torii4}
we extended ${\rm P}_!$ and ${\rm P}^*$
to monoidal presheaf functors.
Let $\mathcal{O}$ be an $\infty$-operad.
We denote by
${\rm Mon}_{\mathcal{O}}^{\rm lax,R}(\cat)$
the $\infty$-category of small
$\mathcal{O}$-monoidal $\infty$-categories
and right adjoint lax $\mathcal{O}$-monoidal functors.
We also denote by
${\rm Mon}_{\mathcal{O}}^{\rm oplax,L}(\cat)$
the $\infty$-category of small $\mathcal{O}$-monoidal
$\infty$-categories and left adjoint
lax $\mathcal{O}$-monoidal functors.
We constructed functors
\[ \begin{array}{llcl}
     \mathbb{P}_{\mathcal{O}!}:&
     {\rm Mon}_{\mathcal{O}}^{\rm lax,R}(\cat)
     & \longrightarrow &
       {\rm Mon}_{\mathcal{O}}^{\rm lax}({\rm Pr}^{\rm L}),\\
     \mathbb{P}_{\mathcal{O}}^*:&
     {\rm Mon}_{\mathcal{O}}^{\rm oplax,L}(\cat)^{\rm op}
     & \longrightarrow &
       {\rm Mon}_{\mathcal{O}}^{\rm lax}({\rm Pr}^{\rm L}),\\
\end{array}\]
where ${\rm Mon}_{\mathcal{O}}^{\rm lax}({\rm Pr}^{\rm L})$
is the $\infty$-category of $\mathcal{O}$-monoidal
$\infty$-categories in the symmetric monoidal
$\infty$-category ${\rm Pr}^{\rm L}$
and lax $\mathcal{O}$-monoidal functors.
For a small $\mathcal{O}$-monoidal
$\infty$-category $\mathcal{C}$,
we denote by
$\mathbb{P}_{\mathcal{O}}(\mathcal{C})$
the $\mathcal{O}$-monoidal $\infty$-category
${\rm Fun}^{\mathcal{O}}(\mathcal{C},\mathcal{S})$
constructed by using Day convolution product
(cf.~\cite[\S2.2.6]{Lurie2}).
We have $\mathbb{P}_{\mathcal{O}!}(\mathcal{C})
=\mathbb{P}(\mathcal{C})
=\mathbb{P}_{\mathcal{O}}^*(\mathcal{C})$.
Let $f: \mathcal{C}\to\mathcal{D}$
be a right adjoint lax
$\mathcal{O}$-monoidal functor,
and let $g: \mathcal{D}\to \mathcal{C}$ 
be a left adjoint oplax
$\mathcal{O}$-monoidal functor.
For simplicity,
we write $f_!$ and $g^*$
for $\mathbb{P}_{\mathcal{O}!}(f)$
and $\mathbb{P}_{\mathcal{O}}^*(g)$,
respectively.

The induced functors on mapping spaces
\[ \begin{array}{llcl}
  \mathbb{P}_{\mathcal{O}!,(\mathcal{C},\mathcal{D})}:&
         {\rm Map}_{{\rm Mon}_{\mathcal{O}}^{\rm lax,R}(\cat)}
         (\mathcal{C},\mathcal{D})
   &\longrightarrow&
         {\rm Map}_{{\rm Mon}_{\mathcal{O}}^{\rm lax}({\rm Pr}^{\rm L})}
         (\mathbb{P}_{\mathcal{O}}(\mathcal{C}),
         \mathbb{P}_{\mathcal{O}}(\mathcal{D})),\\
  \mathbb{P}_{\mathcal{O},(\mathcal{C},\mathcal{D})}^*:&
         {\rm Map}_{{\rm Mon}_{\mathcal{O}}^{\rm oplax,L}(\cat)}
         (\mathcal{C},\mathcal{D})
   &\longrightarrow&
         {\rm Map}_{{\rm Mon}_{\mathcal{O}}^{\rm lax}({\rm Pr}^{\rm L})}
         (\mathbb{P}_{\mathcal{O}}(\mathcal{D}),
         \mathbb{P}_{\mathcal{O}}(\mathcal{C}))\\
\end{array}\]
are fully faithful for any small
$\mathcal{O}$-monoidal $\infty$-categories
$\mathcal{C}$ and $\mathcal{D}$
by \cite[Definitions~2.10 and 2.14]{Torii4}.

For any right adjoint lax
$\mathcal{O}$-monoidal functor $f: \mathcal{C}\to \mathcal{D}$,
we have a left adjoint oplax $\mathcal{O}$-monoidal
functor $f^L: \mathcal{D}\to \mathcal{C}$.
Dually,
for any left adjoint oplax
$\mathcal{O}$-monoidal functor $g: \mathcal{D}\to \mathcal{C}$,
we have a right adjoint lax $\mathcal{O}$-monoidal
functor $g^R: \mathcal{C}\to \mathcal{D}$.
By \cite[Lemma~2.16]{Torii4},
we have $(f^L)^*\simeq f_!$ and $(g^R)_!\simeq g^*$.


\section{Uniqueness of equivalences}
\label{section:uniqueness}

In this section we will give a 
uniqueness result of equivalences
between ${\rm Mon}_{\mathcal{O}}^{\rm lax,R}(\cat)$
and ${\rm Mon}_{\mathcal{O}}^{\rm oplax,L}(\cat)^{\rm op}$
which are compatible with the monoidal
presheaf functors $\mathbb{P}_{\mathcal{O}!}$
and $\mathbb{P}_{\mathcal{O}}^*$.
Using this result,
we will prove that the equivalence of \cite{Torii4}
is equivalent to the restriction
of the equivalence of \cite{HHLN1}
in the next section.

First,
we consider the space
${\rm Map}_{\widehat{\rm Cat}_{\infty/{{\rm Mon}_{\mathcal{O}}^{\rm lax}}
   ({\rm Pr}^{\rm L})}}
(\mathbb{P}_{\mathcal{O}!},
 \mathbb{P}_{\mathcal{O}}^*)$
of functors from 
${\rm Mon}_{\mathcal{O}}^{\rm lax,R}(\cat)$
to ${\rm Mon}_{\mathcal{O}}^{\rm oplax,L}(\cat)^{\rm op}$
which are compatible with 
$\mathbb{P}_{\mathcal{O}!}$
and $\mathbb{P}_{\mathcal{O}}^*$.
By Lemma~\ref{lemma:general-discrete-over-category} below,
we see that it is discrete and
equivalent to the space
${\rm Map}_{\widehat{\rm Cat}_{\infty/{h{\rm Mon}_{\mathcal{O}}^{\rm lax}}({\rm Pr}^{\rm L})}}
(h\mathbb{P}_{\mathcal{O}!},
 h\mathbb{P}_{\mathcal{O}}^*)$
of functors between 
the homotopy categories 
$h{\rm Mon}_{\mathcal{O}}^{\rm lax,R}(\cat)$
and $h{\rm Mon}_{\mathcal{O}}^{\rm oplax,L}(\cat)^{\rm op}$
which are compatible with
$h\mathbb{P}_{\mathcal{O}!}$ and
$h\mathbb{P}_{\mathcal{O}}^*$.

\begin{lemma}\label{lemma:general-discrete-over-category}
Let $P: \mathcal{X}\to \mathcal{Z}$
and $Q:\mathcal{Y}\to\mathcal{Z}$
be functors of $\infty$-categories.
We assume that the induced functor
${\rm Map}_{\mathcal{Y}}(y,y')\to
{\rm Map}_{\mathcal{Z}}(Q(y),Q(y'))$
on mapping spaces
is $(-1)$-truncated
for any objects $y,y'\in\mathcal{Y}$.
Then the space
${\rm Map}_{{\rm Cat}_{\infty/\mathcal{Z}}}
(P,Q)$
is discrete,
and the map
${\rm Map}_{{\rm Cat}_{\infty/\mathcal{Z}}}
(P,Q)
\to
   {\rm Map}_{{\rm Cat}_{\infty/h\mathcal{Z}}}
(hP,hQ)$
is an equivalence.
\end{lemma}

\proof
First,
we shall show
that the space
${\rm Map}_{{\rm Cat}_{\infty/\mathcal{Z}}}
(P,Q)$
is discrete.
For this,
it suffices to show that
the functor $Q$ is a $0$-truncated morphism in $\cat$. 
We consider the diagonal functor
$\Delta: \mathcal{Y}\to \mathcal{Y}
\times_{\mathcal{Z}}\mathcal{Y}$.
By \cite[Lemma~5.5.6.15]{Lurie1},
it suffices to show that
$\Delta$ is $(-1)$-truncated in $\cat$,
that is,
the map $\Delta_*: {\rm Map}_{\cat}(\mathcal{E},\mathcal{Y})
\to {\rm Map}_{\cat}(\mathcal{E},\mathcal{Y}
\times_{\mathcal{Z}}\mathcal{Y})$
is $(-1)$-truncated in $\mathcal{S}$
for any $\infty$-category $\mathcal{E}$.
Note that a map between spaces
is $(-1)$-truncated if and only if
it is fully faithful.
By the assumption,
we see that the functor $\Delta$
is fully faithful.
This implies that
$\Delta_*$ is also fully faithful.

Next,
we shall show that 
the map
${\rm Map}_{{\rm Cat}_{\infty/\mathcal{Z}}}
(P,Q)
\to
   {\rm Map}_{{\rm Cat}_{\infty/h\mathcal{Z}}}
(hP,hQ)$
is an equivalence.
By the assumption,
the following commutative diagram
\[ \begin{array}{ccc}
     \mathcal{Y} &
     \stackrel{\pi_{\mathcal{Y}}}{\longrightarrow} & h\mathcal{Y} \\
    \mbox{$\scriptstyle Q$}
    \bigg\downarrow
    \phantom{\mbox{$\scriptstyle Q$}} 
    & & 
    \phantom{\mbox{$\scriptstyle hQ$}}
    \bigg\downarrow
    \mbox{$\scriptstyle hQ$} \\
    \mathcal{Z} &
    \stackrel{\pi_{\mathcal{Z}}}{\longrightarrow} &
    h\mathcal{Z} \\ 
   \end{array}\]
is pullback.
This induces an equivalence
${\rm Map}_{{\rm Cat}_{\infty/\mathcal{Z}}}(P,Q)
   \stackrel{\simeq}{\to}
   {\rm Map}_{{\rm Cat}_{\infty/h\mathcal{Z}}}
(\pi_{\mathcal{Z}}\circ P,hQ)$.
The desired equivalence follows from 
the fact that 
the functor $\pi_{\mathcal{X}}:
\mathcal{X}\to h\mathcal{X}$
induces an equivalence
${\rm Map}_{{\rm Cat}_{\infty}}(h\mathcal{X},\mathcal{C})
\stackrel{\simeq}{\to}
{\rm Map}_{{\rm Cat}_{\infty}}(\mathcal{X},\mathcal{C})$
for any ordinary category $\mathcal{C}$
by \cite[Proposition~1.2.3.1]{Lurie1}.
\qed

\bigskip

We shall define an object
\[ \mathbf{eq}=({\rm eq},{\rm id})\in 
   {\rm Map}_{\widehat{\rm Cat}_{\infty/{h{\rm Mon}_{\mathcal{O}}^{\rm lax}}({\rm Pr}^{\rm L})}}
   (h\mathbb{P}_{\mathcal{O}!},
    h\mathbb{P}_{\mathcal{O}}^*). \]
Let 
${\rm eq}: h{\rm Mon}_{\mathcal{O}}^{\rm lax,R}(\cat)
\to h{\rm Mon}_{\mathcal{O}}^{\rm oplax, L}(\cat)^{\rm op}$
be a functor which is identity on objects and
assigns to a connected component containing
a right adjoint lax $\mathcal{O}$-monoidal functor $f$
the connected component containing
its left adjoint oplax $\mathcal{O}$-monoidal functor $f^L$.
Since $(f^L)^*\simeq f_!$,
we have the identity natural equivalence
${\rm id}: h\mathbb{P}_{\mathcal{O}!} 
\stackrel{\simeq}{\to}h\mathbb{P}_{\mathcal{O}}^*\circ {\rm eq}$.

Now, 
we will give a condition for
an object of 
${\rm Map}_{\widehat{\rm Cat}_{\infty/{\rm Mon}_{\mathcal{O}}^{\rm lax}({\rm Pr}^{\rm L})}}
(\mathbb{P}_{\mathcal{O}!},\mathbb{P}_{\mathcal{O}}^*)$
to determine an object
of 
${\rm Map}_{\widehat{\rm Cat}_{\infty/h{\rm Mon}_{\mathcal{O}}^{\rm lax}({\rm Pr}^{\rm L})}}
(h\mathbb{P}_{\mathcal{O}!},h\mathbb{P}_{\mathcal{O}}^*)$
which is equivalent to $\mathbf{eq}$.

\begin{lemma}\label{lemma:characterizetion-eq}
Let $\mathbf{F}=(F,\delta)$ be an object
of ${\rm Map}_{\widehat{\rm Cat}_{\infty/{\rm Mon}_{\mathcal{O}}^{\rm lax}({\rm Pr}^{\rm L})}}
(\mathbb{P}_{\mathcal{O}!},\mathbb{P}_{\mathcal{O}}^*)$,
where $F$ is a functor from ${\rm Mon}_{\mathcal{O}}^{\rm lax,R}(\cat)$
to ${\rm Mon}_{\mathcal{O}}^{\rm oplax,L}(\cat)^{\rm op}$
and $\delta: \mathbb{P}_{\mathcal{O}!}
\stackrel{\simeq}{\to}
\mathbb{P}_{\mathcal{O}}^*\circ F$
is a natural equivalence. 
Then $\mathbf{F}$ determines
an object of 
${\rm Map}_{\widehat{\rm Cat}_{\infty/h{\rm Mon}_{\mathcal{O}}^{\rm lax}({\rm Pr}^{\rm L})}}
(h\mathbb{P}_{\mathcal{O}!},h\mathbb{P}_{\mathcal{O}}^*)$
which is equivalent to $\mathbf{eq}$
if and only if
there is an equivalence $\zeta_{\mathcal{C}}: 
F(\mathcal{C})\to \mathcal{C}$ 
for each small $\mathcal{O}$-monoidal $\infty$-category $\mathcal{C}$
such that $\zeta_{\mathcal{C}}^*:
\mathbb{P}_{\mathcal{O}}(\mathcal{C})\to 
\mathbb{P}_{\mathcal{O}}(F(\mathcal{C}))$
is equivalent to $\delta_{\mathcal{C}}$.
\end{lemma}

\proof
First, 
we suppose that 
$h\mathbf{F}=(hF,h\delta)$
is equivalent to $\mathbf{eq}$.
Then there is a natural equivalence
$\tau: hF\stackrel{\simeq}{\to}{\rm eq}$
which makes the following diagram commute
\begin{equation}\label{eq:id-h-delta-PO-tau-triangle}
\xymatrix{
    & h\mathbb{P}_{\mathcal{O}!}
      \ar[dl]_{{\rm id}} \ar[dr]^{h\delta} & \\
    h\mathbb{P}_{\mathcal{O}}^*\circ {\rm eq}
    \ar[rr]^{h\mathbb{P}_{\mathcal{O}}^*\circ \tau}
    &&
    h\mathbb{P}_{\mathcal{O}}^*\circ hF.\\
}\end{equation}
For each $\mathcal{C}$,
we obtain a desired equivalence $\zeta_{\mathcal{C}}$
by taking a representative of the component $\tau_{\mathcal{C}}$.

Conversely,
we suppose that we have
an equivalence $\zeta_{\mathcal{C}}:F(\mathcal{C})
\stackrel{\simeq}{\to}\mathcal{C}$
such that $\zeta_{\mathcal{C}}^*\simeq\delta_{\mathcal{C}}$
for each $\mathcal{C}$.
We will show that 
the collection $\{\zeta_{\mathcal{C}}\}$ of equivalences
determines a natural 
equivalence $\tau: hF\stackrel{\simeq}{\to}{\rm eq}$. 
We consider a diagram $\sigma$:
\[   \begin{array}{ccc}
    F(\mathcal{C}) 
    & \stackrel{\zeta_{\mathcal{C}}}{\to} & 
    \mathcal{C}\\
    \mbox{$\scriptstyle F(f)$}
    \bigg\uparrow
    \phantom{\mbox{$\scriptstyle F(f)$}} 
    & & 
    \phantom{\mbox{$\scriptstyle f^L$}}
    \bigg\uparrow
    \mbox{$\scriptstyle f^L$} \\
    F(\mathcal{D})    
    & \stackrel{\zeta_{\mathcal{D}}}{\to } & 
    \mathcal{D} \\
   \end{array}
\]
in ${\rm Mon}_{\mathcal{O}}^{\rm oplax,L}(\cat)$
for any morphism $f:\mathcal{C}\to \mathcal{D}$
in ${\rm Mon}_{\mathcal{O}}^{\rm lax,R}(\cat)$.
By the assumption that
$\zeta_{\mathcal{C}}^*\simeq \delta_{\mathcal{C}}$
and the fact that $(f^L)^*\simeq f_!$,
we see that the diagram
$\mathbb{P}_{\mathcal{O}}^*(\sigma)$
is homotopy commutative.
Since 
$\mathbb{P}_{\mathcal{O},(F(\mathcal{D}),\mathcal{C})}^*$
is fully faithful,
we see that the diagram $\sigma$
is also homotopy commutative. 
Hence the collection $\{\zeta_{\mathcal{C}}\}$ determines
a natural equivalence
$\tau: hF\stackrel{\simeq}{\to} {\rm eq}$.
By the construction,
we obtain a commutative diagram~(\ref{eq:id-h-delta-PO-tau-triangle}).
Therefore,
$h\mathbf{F}=(hF,h\delta)$ is equivalent to 
$\mathbf{eq}=({\rm eq},{\rm id})$.
\qed

\bigskip

By \cite[Theorem~4.4]{Torii4},
we have an object $(T,\theta)$
of ${\rm Map}_{\widehat{\rm Cat}_{\infty/{\rm Mon}_{\mathcal{O}}^{\rm lax}({\rm Pr}^{\rm L})}}
(\mathbb{P}_{\mathcal{O}!},\mathbb{P}_{\mathcal{O}}^*)$.
We shall show that 
$(T,\theta)$
induces an object equivalent to $\mathbf{eq}$
in ${\rm Map}_{\widehat{\rm Cat}_{\infty/
{h\rm Mon}_{\mathcal{O}}^{\rm lax}({\rm Pr}^{\rm L})}}
(h\mathbb{P}_{\mathcal{O}!},h\mathbb{P}_{\mathcal{O}}^*)$.

\begin{proposition}\label{prop:T-determine-eq}
The object $(T,\theta)$
determines an object 
of ${\rm Map}_{\widehat{\rm Cat}_{\infty/h{\rm Mon}_{\mathcal{O}}^{\rm lax}({\rm Pr}^{\rm L})}}
(h\mathbb{P}_{\mathcal{O}!},h\mathbb{P}_{\mathcal{O}}^*)$
which is equivalent to $\mathbf{eq}$.
\end{proposition}

\proof
We use the notation in \cite[\S4]{Torii4}.
For a small $\mathcal{O}$-monoidal 
$\infty$-category $\mathcal{D}$,
we take  
a left and right universal object 
$(\mathcal{C},\mathcal{D},f)$ 
of $\mathcal{M}_{\mathcal{O}}^{\rm lax,R}$
with $f\simeq g_!$,
where $g:\mathcal{C}\stackrel{\simeq}{\to} \mathcal{D}$
is an equivalence in 
${\rm Mon}_{\mathcal{O}}^{\rm lax,R}(\cat)$.
Then we have
an equivalence $\mu:
T(\mathcal{D})\stackrel{\simeq}{\to}\mathcal{C}$
in ${\rm Mon}_{\mathcal{O}}^{\rm oplax,L}(\cat)$
which makes the following diagram commute
\[ \xymatrix{
    \mathbb{P}_{\mathcal{O}}(\mathcal{C})
    \ar[rr]^{\mu^*}\ar[dr]_{f}
    &&
    \mathbb{P}_{\mathcal{O}}(T(\mathcal{D}))
    \\
    &\mathbb{P}_{\mathcal{O}}(\mathcal{D}).
     \ar[ur]_{\theta_{\mathcal{D}}}& \\
}\]
We set $\zeta_{\mathcal{D}}=g\circ \mu: 
T(\mathcal{D})\to \mathcal{D}$.
Then $\zeta_{\mathcal{D}}$ is an equivalence
and $\zeta_{\mathcal{D}}^*\simeq \theta_{\mathcal{D}}$ 
since $g^*\simeq f^{-1}$.
The desired result follows from
Lemma~\ref{lemma:characterizetion-eq}.
\qed

\bigskip

We define a space ${\mathcal Eq}$ by
\[  {\mathcal Eq} 
  = {\rm Map}_{\widehat{\rm Cat}_{\infty/
    {\rm Mon}_{\mathcal{O}}^{\rm lax}({\rm Pr}^{\rm L})}}
   (\mathbb{P}_{\mathcal{O}!},\mathbb{P}_{\mathcal{O}}^*)
   \times_{
   {\rm Map}_{\widehat{\rm Cat}_{\infty/
    h{\rm Mon}_{\mathcal{O}}^{\rm lax}({\rm Pr}^{\rm L})}}
    (h\mathbb{P}_{\mathcal{O}!},h\mathbb{P}_{\mathcal{O}}^*)}
   \{
     \mathbf{eq}
    \}.
\]
We
can regard ${\mathcal Eq}$
as the space of equivalences
between ${\rm Mon}_{\mathcal{O}}^{\rm lax,R}(\cat)$
and ${\rm Mon}_{\mathcal{O}}^{\rm oplax,L}(\cat)^{\rm op}$
which are compatible with 
the monoidal presheaf functors
$\mathbb{P}_{\mathcal{O}!}$ and $\mathbb{P}_{\mathcal{O}}^*$.
By Lemma~\ref{lemma:general-discrete-over-category},
we obtain the following theorem. 

\begin{theorem}\label{thm:uniqueness-equivalences}
The space ${\mathcal Eq}$ is contractible.
\end{theorem}

\section{Canonical equivalence 
between $(T,\theta)$ and $(H,\gamma)$}
\label{section:canonical-equivalence}

In this section
we shall show that 
the equivalence $T$ in \cite{Torii4} 
is equivalent to 
the restriction $H$ of the equivalence in \cite{HHLN1}.
The functor $T$ equipped with
the commutative diagram in \cite[Theorem~4.4]{Torii4}
gives rise to an object $(T,\theta)$
in ${\mathcal Eq}$
by Proposition~\ref{prop:T-determine-eq}.
Hence it suffices to show that
$H$ extends to an object $(H,\gamma)$ of ${\mathcal Eq}$
by Theorem~\ref{thm:uniqueness-equivalences}. 

First,
we have to show that 
$H$ gives a morphism from $\mathbb{P}_{\mathcal{O}!}$
to $\mathbb{P}_{\mathcal{O}}^*$ in
$\widehat{\rm Cat}_{\infty/{\rm Mon}_{\mathcal{O}}^{\rm lax}({\rm Pr}^{\rm L})}$.
In order to prove this,
we will show that there are natural transformations
$I\to \mathbb{P}_{\mathcal{O}!}$
and 
$I\to \mathbb{P}_{\mathcal{O}}^*\circ H$
of functors from
${\rm Mon}_{\mathcal{O}}^{\rm lax,R}(\cat)$
to ${\rm Mon}_{\mathcal{O}}^{\rm lax}(\wcat)$
whose components are equivalent to Yoneda embeddings,
where $I$ is the inclusion functor
${\rm Mon}_{\mathcal{O}}^{\rm lax,R}(\cat)\to 
{\rm Mon}_{\mathcal{O}}^{\rm lax}(\wcat)$.
After that,
we will show that these natural transformations
are canonically equivalent.

We consider a pullback diagram
\[ \begin{array}{ccc}
   \mathcal{F}& \longrightarrow &
   {\rm Fun}({\rm Mon}_{\mathcal{O}}^{\rm lax,R}(\cat),
             {\rm Mon}_{\mathcal{O}}^{\rm lax}({\rm Pr}^{\rm L}))\\
   \bigg\downarrow & & \bigg\downarrow \\
   {\rm Fun}({\rm Mon}_{\mathcal{O}}^{\rm lax,R}(\cat),
  {\rm Mon}_{\mathcal{O}}^{\rm lax}(\wcat)
)_{I/}
    & \longrightarrow &
   {\rm Fun}({\rm Mon}_{\mathcal{O}}^{\rm lax,R}(\cat),
  {\rm Mon}_{\mathcal{O}}^{\rm lax}(\wcat)
). \\
   \end{array}     \]
We need the following lemma which  
characterizes initial objects of $\mathcal{F}$.

\begin{lemma}\label{lemma:characterization-initial-F}
An object $(\sigma: I\to F)$ of $\mathcal{F}$
is initial if
the component $\sigma_{\mathcal{C}}:
\mathcal{C}\to F(\mathcal{C})$
is equivalent to the Yoneda embedding
$J_{\mathcal{C}}:\mathcal{C}\to\mathbb{P}_{\mathcal{O}}(\mathcal{C})$
for each $\mathcal{C}\in 
{\rm Mon}_{\mathcal{O}}^{\rm lax,R}(\cat)$,
that is,
there is an equivalence
$F(\mathcal{C})\stackrel{\simeq}{\to}
\mathbb{P}_{\mathcal{O}}(\mathcal{C})$
which makes the following diagram commute
\[ \xymatrix{
     &\mathcal{C}
      \ar[dl]_{\sigma_{\mathcal{C}}}\ar[dr]^{J_{\mathcal{C}}} & \\
      F(\mathcal{C})\ar[rr]&&
      \mathbb{P}_{\mathcal{O}}(\mathcal{C}).
}\]
\end{lemma} 

\proof
For any object $(\tau: I\to G)$ of $\mathcal{F}$,
we have the following pullback diagrams
\[ \begin{array}{ccc}
   {\rm Map}_{\mathcal{F}}(\sigma,\tau)
   &\longrightarrow&
   {\rm Map}_{{\rm Fun}({\rm Mon}_{\mathcal{O}}^{\rm lax,R}(\cat),
               {\rm Mon}_{\mathcal{O}}^{\rm lax}({\rm Pr}^{\rm L}))}
   (\sigma,\tau)    \\[2mm]
   \bigg\downarrow & & \bigg\downarrow \\
   {\rm Map}_{{\rm Fun}({\rm Mon}_{\mathcal{O}}^{\rm lax,R}(\cat),
     {\rm Mon}_{\mathcal{O}}^{\rm lax}(
     \wcat
     ))_{I/}}   
   (\sigma,\tau) 
   & \longrightarrow &
   {\rm Map}_{{\rm Fun}({\rm Mon}_{\mathcal{O}}^{\rm lax,R}(\cat),
     {\rm Mon}_{\mathcal{O}}^{\rm lax}(
     \wcat
     ))}   
   (\sigma,\tau)    \\
   \bigg\downarrow & & \bigg\downarrow \\
   * & \longrightarrow &
   {\rm Map}_{{\rm Fun}({\rm Mon}_{\mathcal{O}}^{\rm lax,R}(\cat),
     {\rm Mon}_{\mathcal{O}}^{\rm lax}(
      \wcat
     ))}   
   (I,\tau).    \\
   \end{array}\]
In order to prove that the composite of the left vertical arrows
is an equivalence,
it suffices to show that
the composite of the right vertical arrows is an equivalence.
The composite of the right vertical arrows is identified with
a map between ends
\[ 
   \int_{\mathcal{C}\in {\rm Mon}_{\mathcal{O}}^{\rm lax,R}(\cat)}
   {\rm Map}_{{\rm Mon}_{\mathcal{O}}^{\rm lax}({\rm Pr}^{\rm L})}
            (F(\mathcal{C}),G(\mathcal{C})) 
   \longrightarrow 
   \int_{\mathcal{C}\in {\rm Mon}_{\mathcal{O}}^{\rm lax,R}(\cat)}
   {\rm Map}_{{\rm Mon}_{\mathcal{O}}^{\rm lax}(\wcat)}
            (\mathcal{C},G(\mathcal{C})),  \]
which is an equivalence
by the assumption and
\cite[Proposition~2.7]{Torii4}.
\qed

\bigskip

Next, we shall show that there is an
object $(\alpha: I\to \mathbb{P}_{\mathcal{O}!})$ of $\mathcal{F}$
whose components are equivalent to Yoneda embeddings.

\begin{lemma}\label{lemma:PO-extension-to-I}
The functor $\mathbb{P}_{\mathcal{O}!}$
extends to an object
$(\alpha: I\to \mathbb{P}_{\mathcal{O}_!})$
in $\mathcal{F}$
such that 
$\alpha_{\mathcal{C}}\simeq J_{\mathcal{C}}$
for each $\mathcal{C}\in
{\rm Mon}_{\mathcal{O}}^{\rm lax,R}(\cat)$.
\end{lemma}

\proof
Let $i: {\rm Mon}_{\mathcal{O}}^{\rm lax}({\rm Pr}^{\rm L})
\to {\rm Mon}_{\mathcal{O}}^{\rm lax}(\wcat)$
be the inclusion functor.
We have an adjunction
\[ i_!: {\rm Fun}
        ({\rm Mon}_{\mathcal{O}}^{\rm lax}({\rm Pr}^{\rm L}),
         \widetilde{S})
        \leftrightarrows
          {\rm Fun}
        ({\rm Mon}_{\mathcal{O}}^{\rm lax}(\wcat),
         \widetilde{S}) : i^*,\]
where $i^*$ is the restriction functor and
$i_!$ is a left Kan extension along $i$.
The counit of the adjunction
$(i_!i^*\to {\rm id})$
gives a functor
${\rm Fun}({\rm Mon}_{\mathcal{O}}^{\rm lax}(\wcat),\widetilde{S})
\times [1]
   \to
   {\rm Fun}({\rm Mon}_{\mathcal{O}}^{\rm lax}(\wcat),\widetilde{S})$.
Composing with 
${\rm Mon}_{\mathcal{O}}^{\rm lax,R}(\cat)^{\rm op}
\to
  {\rm Fun}({\rm Mon}_{\mathcal{O}}^{\rm lax}(\wcat),\widetilde{S})$
which assigns 
${\rm Map}_{{\rm Mon}_{\mathcal{O}}^{\rm lax}(\wcat)}(\mathcal{C},-)$
to $\mathcal{C}$,
we obtain a functor
${\rm Mon}_{\mathcal{O}}^{\rm lax,R}(\cat)^{\rm op}\times [1]
\to
   {\rm Fun}({\rm Mon}_{\mathcal{O}}^{\rm lax}(\wcat),\widetilde{S})$.
We will show that 
this functor factors through
${\rm Mon}_{\mathcal{O}}^{\rm lax}(\wcat)^{\rm op}$.
For this,
it suffices to show that
$i_!i^*M(\mathcal{C})$
is representable 
for any $\mathcal{C}\in {\rm Mon}_{\mathcal{O}}^{\rm lax,R}(\cat)$,
where $M(\mathcal{C})=
{\rm Map}_{{\rm Mon}_{\mathcal{O}}^{\rm lax}(\wcat)}(\mathcal{C},-)$.
The functor $i^*M(\mathcal{C})$
is equivalent to 
${\rm Map}_{{\rm Mon}_{\mathcal{O}}^{\rm lax}({\rm Pr}^{\rm L})}
(\mathbb{P}_{\mathcal{O}}(\mathcal{C}),-)$
by \cite[Proposition~2.7]{Torii4}.
Hence $i_!i^*M(\mathcal{C})$ is represented by
$\mathbb{P}_{\mathcal{O}}(\mathcal{C})$
since ${\rm Map}_{{\rm Mon}_{\mathcal{O}}^{\rm lax}(\wcat)}
(\mathbb{P}_{\mathcal{O}}(\mathcal{C}),-)$
is a left Kan extension of
${\rm Map}_{{\rm Mon}_{\mathcal{O}}^{\rm lax}({\rm Pr}^{\rm L})}
(\mathbb{P}_{\mathcal{O}}(\mathcal{C}),-)$
along $i$.
Therefore,
we obtain a functor
${\rm Mon}_{\mathcal{O}}^{\rm lax, R}(\cat)\times [1]\to
{\rm Mon}_{\mathcal{O}}^{\rm lax}(\wcat)$,
which is a natural transformation
$\alpha: I \to \mathbb{P}_{\mathcal{O}!}$
such that $\alpha_{\mathcal{C}}\simeq J_{\mathcal{C}}$
for each $\mathcal{C}$.
\qed

\bigskip

Next, we would like to
show that there is an object
$(\beta: I\to \mathbb{P}_{\mathcal{O}}^*\circ H)$
of $\mathcal{F}$
whose components are equivalent to Yoneda embeddings.
For this purpose,
we need some preliminary constructions.

Let $p=(p_A,p_{\mathcal{O}}): Y\to A\times\mathcal{O}^{\otimes}$
be the curved orthofibration associated to 
a functor ${\rm St}(p):
A^{\rm op}\to {\rm Mon}_{\mathcal{O}}^{\rm lax,R}(\cat)$
(see \cite{HHLN1} for the definitions of curved orthofibrations
and (op-)Gray fibrations).
Although 
we will take $p$ as the identity map later, 
it is notationally convenient to introduce
a placeholder.
Note that 
$(p_{\mathcal{O}},p_A):
Y\to \mathcal{O}^{\otimes}\times A$
is a Gray fibration.
We denote by $p^{\wedge}:
Y^{\wedge}\to \mathcal{O}^{\otimes,{\rm op}}\times A$
the op-Gray fibration
obtained from $p$ by coCartesian straightening followed by
Cartesian unstraightening
with respect to $p_{\mathcal{O}}: Y\to \mathcal{O}^{\otimes}$.
Then $p^{\wedge}$ 
is also a curved orthofibration
associated to the functor
$H\circ {\rm St}(p):
A^{\rm op}\to {\rm Mon}_{\mathcal{O}}^{\rm oplax,L}(\cat)^{\rm op}$.
By taking opposite categories,
we obtain a curved orthofibration
$p^{\vee}=(p^{\vee}_A,p^{\vee}_{\mathcal{O}})=(p^{\wedge})^{\rm op}:
Y^{\vee}=(Y^{\wedge})^{\rm op}\to A^{\rm op}\times \mathcal{O}^{\otimes}$
which is associated to
the functor $(-)^{\vee}\circ H\circ {\rm St}(p):
A^{\rm op}\to {\rm Mon}_{\mathcal{O}}^{\rm lax,R}(\cat)^{\rm op}$.

Let ${\rm Tw}^l(A)$ be the twisted arrow $\infty$-category
equipped with a left fibration
${\rm Tw}^l(A)\to A^{\rm op}\times A$.
We consider the following pullback diagram
\[ \begin{array}{ccc}
    E(Y/A) & \longrightarrow & Y^{\vee}
    \times_{\mathcal{O}^{\otimes}} Y \\
    \mbox{$\scriptstyle q$}\bigg\downarrow
    \phantom{\mbox{$\scriptstyle q$}} 
     & & 
    \phantom{\mbox{$\scriptstyle p_A^{\vee}\times p_A$}}
    \bigg\downarrow
    \mbox{$\scriptstyle p_A^{\vee}\times p_A$} \\[3mm]
    {\rm Tw}^l(A) & \longrightarrow &A^{\rm op}\times A. \\
   \end{array} \]
Note that the vertical arrows are Cartesian fibrations,
and that the fiber $E(Y/A)_{\phi}$ is an $\mathcal{O}$-monoidal
$\infty$-category $Y^{\vee}_s\times_{\mathcal{O}^{\otimes}}Y_t$
for any $(\phi: s\to t)\in {\rm Tw}^l(A)$.

Let ${\rm Nat}(I,\mathbb{P}_{\mathcal{O}}^*\circ H)_A$
be 
the space
of natural transformations
of functors from 
$I\circ {\rm St}(p)$
to $\mathbb{P}_{\mathcal{O}}^*\circ H\circ {\rm St}(p)$.
We will show that there is a fully faithful functor
${\rm Nat}(I,\mathbb{P}_{\mathcal{O}}^*\circ H)_A
\to 
{\rm Map}_{\wcat}
(E(Y/A),\mathcal{S})$.

For this purpose,
we will first show that 
the Cartesian symmetric monoidal structure on
the $\infty$-category ${\rm Mon}_{\mathcal{O}}^{\rm lax}(\cat)$
is closed.
By naturality of 
the construction of Day convolution product
\cite[Construction~2.2.6.18]{Lurie2},
we have a functor
\[ {\rm Fun}^{\mathcal{O}}(-,-):
   {\rm Mon}_{\mathcal{O}}^{\rm lax}(\cat)^{\rm op}
   \times
   {\rm Mon}_{\mathcal{O}}^{\rm lax}(\cat)
   \longrightarrow
   {\rm Mon}_{\mathcal{O}}^{\rm lax}(\cat). \]

\begin{lemma}\label{lemma:natural-eq-fun-O}
There is a natural equivalence
\[ {\rm Fun}^{\mathcal{O}}(\mathcal{A},
   {\rm Fun}^{\mathcal{O}}(\mathcal{B},\mathcal{C}))
   \simeq
   {\rm Fun}^{\mathcal{O}}(\mathcal{A}\times_{\mathcal{O}}\mathcal{B},
   \mathcal{C}) \]
of $\mathcal{O}$-monoidal $\infty$-categories
for any $\mathcal{A},\mathcal{B},
\mathcal{C}\in {\rm Mon}_{\mathcal{O}}^{\rm lax}(\cat)$.
This induces a natural equivalence
\[ {\rm Fun}_{\mathcal{O}}^{\rm lax}
   (\mathcal{A},
   {\rm Fun}^{\mathcal{O}}(\mathcal{B},\mathcal{C}))
   \simeq
   {\rm Fun}_{\mathcal{O}}^{\rm lax}
   (\mathcal{A}\times_{\mathcal{O}}\mathcal{B},
   \mathcal{C}) \]
of $\infty$-categories.
\end{lemma}

\proof
Since we have a natural equivalence
${\rm Fun}_{\mathcal{O}}^{\rm lax}(\mathcal{O},
{\rm Fun}^{\mathcal{O}}(\mathcal{X},\mathcal{Y}))
\simeq
{\rm Fun}_{\mathcal{O}}^{\rm lax}(\mathcal{X},\mathcal{Y})$
for any $\mathcal{O}$-monoidal $\infty$-categories
$\mathcal{X}$ and $\mathcal{Y}$
by \cite[Example~2.2.6.3]{Lurie2},
the second part follows from the first part.
Thus,
it suffices to show the first part.

Let 
$L: {\rm Mon}_{\mathcal{O}}^{\rm lax}(\cat)^{\rm op}
      \times
      {\rm Mon}_{\mathcal{O}}^{\rm lax}(\cat)^{\rm op}
      \times
      {\rm Mon}_{\mathcal{O}}^{\rm lax}(\cat)
      \to
      {\rm Mon}_{\mathcal{O}}^{\rm lax}(\cat)$
be the functor given by
$L(\mathcal{A},\mathcal{B},\mathcal{C})
={\rm Fun}^{\mathcal{O}}(A\times_{\mathcal{O}}\mathcal{B},
\mathcal{C})$,
and let 
$R: {\rm Mon}_{\mathcal{O}}^{\rm lax}(\cat)^{\rm op}
      \times
      {\rm Mon}_{\mathcal{O}}^{\rm lax}(\cat)^{\rm op}
      \times
      {\rm Mon}_{\mathcal{O}}^{\rm lax}(\cat)
      \to
      {\rm Mon}_{\mathcal{O}}^{\rm lax}(\cat)$
be the functor given by
$R(\mathcal{A},\mathcal{B},\mathcal{C})
={\rm Fun}^{\mathcal{O}}(A,
{\rm Fun}^{\mathcal{O}}(\mathcal{B},\mathcal{C}))$.
We will compare $L$ and $R$ with a functor
$G: {\rm Mon}_{\mathcal{O}}^{\rm lax}(\cat)^{\rm op}
      \times
      {\rm Mon}_{\mathcal{O}}^{\rm lax}(\cat)^{\rm op}
      \times
      {\rm Mon}_{\mathcal{O}}^{\rm lax}(\cat)
\to
      {\rm Fun}({\rm Mon}_{\mathcal{O}}^{\rm lax}(\cat)^{\rm op},
      \widehat{\mathcal{S}})$
given by
$G(\mathcal{A},\mathcal{B},\mathcal{C})
   =
   (\mathcal{Z}\mapsto {\rm Map}_{{\rm Mon}_{\mathcal{O}}^{\rm lax}(\cat)}
    (\mathcal{Z}\times_{\mathcal{O}}\mathcal{A}
    \times_{\mathcal{O}}\mathcal{B},\mathcal{C}))$.
By the universal property 
of ${\rm Fun}^{\mathcal{O}}(-,-)$,
we see that $j\circ L$ and $j\circ R$
are equivalent to $G$,
where $j: {\rm Mon}_{\mathcal{O}}^{\rm lax}(\cat)
\to {\rm Fun}({\rm Mon}_{\mathcal{O}}^{\rm lax}(\cat)^{\rm op},
\widehat{S})$
is the Yoneda embedding.
Since $j$ is fully faithful,
we obtain that $L$ is equivalent to $R$.
\qed

%

\begin{proposition}
\label{prop:fully-faithful-Nat-Fun}
There is a fully faithful functor
${\rm Nat}(I,\mathbb{P}_{\mathcal{O}}^*\circ H)_A\to
{\rm Map}_{\wcat}
(E(Y/A),\mathcal{S})$
whose essential image is spanned by
those functors $F: E(Y/A)\to \mathcal{S}$
which carry $q$-Cartesian morphisms to
equivalences and whose restriction
$F_{\phi}: E(Y/A)_{\phi}\to \mathcal{S}$
is an $E(Y/A)_{\phi}$-monoid object for 
each $\phi\in {\rm Tw}^l(A)$.
\end{proposition}

\proof
By \cite[Proposition~5.1]{GHN}
(or \cite[Proposition~2.3]{Glasman}),
the 
space
${\rm Nat}(I,\mathbb{P}_{\mathcal{O}}^*\circ H)_A$
is equivalent to the end
$\int_{t\in A}
{\rm Map}_{{\rm Mon}_{\mathcal{O}}^{\rm lax}(\cat)}
  (Y_t,\mathbb{P}_{\mathcal{O}}(Y_t))$.
By using Lemma~\ref{lemma:natural-eq-fun-O},
we have 
natural
equivalences
${\rm Fun}_{\mathcal{O}}^{\rm lax}
(Y_t,\mathbb{P}_{\mathcal{O}}(Y_s))\simeq
{\rm Fun}_{\mathcal{O}}^{\rm lax}
(Y_s^{\vee}\times_{\mathcal{O}^{\otimes}}Y_t,
\mathcal{S}_{\mathcal{O}}^{\otimes})
\simeq
{\rm Fun}^{\times}
(Y_s^{\vee}\times_{\mathcal{O}^{\otimes}}Y_t,
\mathcal{S})$,
where ${\rm Fun}^{\times}(Y_s^{\vee}\times_{\mathcal{O}^{\otimes}}Y_t,
\mathcal{S})$
is a full subcategory of
${\rm Fun}(Y_s^{\vee}\times_{\mathcal{O}^{\otimes}}Y_t,
\mathcal{S})$
spanned by $Y_s^{\vee}\times_{\mathcal{O}^{\otimes}}Y_t$-monoid objects
in $\mathcal{S}$.
Hence we obtain 
a fully faithful functor
${\rm Nat}(I,\mathbb{P}_{\mathcal{O}}^*\circ H)_A\to
\int_{t\in A}
{\rm Map}_{\wcat}
(Y_t^{\vee}\times_{\mathcal{O}^{\otimes}}Y_t,\mathcal{S})
\simeq
{\rm Map}_{\wcat}
(\int^{t\in A}Y_t^{\vee}\times_{\mathcal{O}^{\otimes}}Y_t,
\mathcal{S})$,
where $\int^{t\in A}Y_t^{\vee}\times_{\mathcal{O}^{\otimes}}Y_t$
is a coend in $\cat$.

The coend $\int^{t\in A}(Y^{\vee})_t\times_{\mathcal{O}^{\otimes}}Y_t$
is a colimit of the composite functor
${\rm Tw}^l(A)\to A^{\rm op}\times A
\stackrel{\chi}{\longrightarrow} \cat$,
where $\chi$ classifies the Cartesian fibration
$p_A^{\vee}\times p_A:
Y^{\vee}\times_{\mathcal{O}^{\otimes}}Y\to A^{\rm op}\times A$. 
Hence it is obtained from $E(Y/A)$
by inverting
all $q$-Cartesian morphisms
by \cite[Corollary~3.3.4.3]{Lurie1}.
Therefore,
there is a fully faithful functor
$
{\rm Map}_{\wcat}
(\int^{t\in A}Y_t^{\vee}\times_{\mathcal{O}^{\otimes}}Y_t,
\mathcal{S})\to
{\rm Map}_{\wcat}
(E(Y/A),\mathcal{S})$
whose essential image is spanned by
those functors which carry $q$-Cartesian morphisms
to equivalences.
By composing these two fully faithful functors,
we obtain the desired functor.
\qed

\bigskip

Next, 
we will construct a functor $B: E(Y/A)\to \mathcal{S}$
which lies in the essential image
of the functor
${\rm Nat}(I,\mathbb{P}_{\mathcal{O}}^*\circ H)_A
\to 
{\rm Map}_{\wcat}
(E(Y/A),\mathcal{S})$
in 
Proposition~\ref{prop:fully-faithful-Nat-Fun}.

By \cite[\S5]{BGN},
we have a fiberwise mapping space functor
${\rm Map}_{Y/\mathcal{O}^{\otimes}}:
   Y^{\vee}\times_{\mathcal{O}^{\otimes}}Y
\to
\mathcal{S}$
and an associated left fibration
$\widetilde{O}(Y/\mathcal{O}^{\otimes})
\to
   Y^{\vee}\times_{\mathcal{O}^{\otimes}}Y$.
By the naturality of construction,
we have a map
$\widetilde{O}(Y/\mathcal{O}^{\otimes})
\to
   \widetilde{O}(A/*)\simeq {\rm Tw}^l(A)$
which makes the following diagram commute
\[ \begin{array}{ccc}
    \widetilde{O}(Y/\mathcal{O}^{\otimes})
    &\longrightarrow & Y^{\vee}\times_{\mathcal{O}^{\otimes}}Y\\[1mm]
   \bigg\downarrow & & 
   \phantom{\mbox{$\scriptstyle p_A^{\vee}\times p_A$}}
   \bigg\downarrow
   \mbox{$\scriptstyle p_A^{\vee}\times p_A$} \\[3mm]
   {\rm Tw}^l(A) & \longrightarrow & A^{\rm op}\times A.\\
   \end{array}\]
Hence we obtain a functor
$\widetilde{O}(Y/\mathcal{O}^{\otimes}) 
\to
   E(Y/A)$.

We will show that 
$\widetilde{O}(Y/\mathcal{O}^{\otimes}) \to
   E(Y/A)$ is a left fibration.
The functor $E(Y/A)\to Y^{\vee}\times_{\mathcal{O}^{\otimes}}Y$
is a left fibration
since it is obtained by pullback from
the left fibration ${\rm Tw}^l(A)\to A^{\rm op}\times A$. 
Since the composite
$\widetilde{O}(Y/\mathcal{O}^{\otimes})
\to E(Y/A)\to Y^{\vee}\times_{\mathcal{O}^{\otimes}}Y$
is a left fibration,
the functor
$\widetilde{O}(Y/\mathcal{O}^{\otimes})
\to E(Y/A)$ is also a left fibration.
Hence we obtain an associated functor
$B: E(Y/A)\to \mathcal{S}$
which assigns to $((y,y'),\phi)$
with $(y,y')\in Y^{\vee}\times_{\mathcal{O}^{\otimes}} Y$
and
$\phi: p_A^{\vee}(y)\to p_A(y')$
the mapping space
\[ {\rm Map}_{Y_o}^{\phi}(y,y')=
   {\rm Map}_{Y_o}(y,y')\times_{{\rm Map}_{A}(p_A^{\vee}(y),p_A(y'))}
   \{\phi\} ,\]
where $o=p_{\mathcal{O}}^{\vee}(y)=p_{\mathcal{O}}(y')\in
\mathcal{O}^{\otimes}$.

We will show that the functor 
$B: E(Y/A)\to \mathcal{S}$
lies in the essential image of 
the fully faithful functor
${\rm Nat}(I,\mathbb{P}_{\mathcal{O}}^*\circ H)_A
\to
{\rm Map}_{\wcat}
(E(Y/A),\mathcal{S})$
in 
Proposition~\ref{prop:fully-faithful-Nat-Fun}.

\begin{lemma}\label{lemma:PO-star-extesion}
The functor $B: E(Y/A)\to \mathcal{S}$
carries $q$-Cartesian morphisms to equivalences.
\end{lemma}

\proof
We take a $q$-Cartesian morphism 
$\Phi:
((y,y'),\phi)\to ((\widetilde{y},\widetilde{y}'),\widetilde{\phi})$
of $E(Y/A)$.
We may assume that
$o=p_{\mathcal{O}}^{\vee}(y)=p_{\mathcal{O}}(y')
  =p_{\mathcal{O}}^{\vee}(\widetilde{y})=p_{\mathcal{O}}(\widetilde{y}')
\in\mathcal{O}^{\otimes}$.
We will show that $\Phi$
induces an equivalence between
${\rm Map}_{Y_o}^{\phi}(y,y')$
and 
${\rm Map}_{Y_o}^{\widetilde{\phi}}(\widetilde{y},\widetilde{y}')$.

We denote by $p_o :Y_o\to A$ the restriction of $p$ 
at $o\in\mathcal{O}^{\otimes}$
which is a Cartesian and coCartesian fibration.
The morphism $\Phi$ consists of
a $p_o$-coCartesian morphism $f: \widetilde{y}\to y$,
a $p_o$-Cartesian morphism
$g: y'\to \widetilde{y}'$,
and a commutative diagram 
\[ \begin{array}{ccc}
    p_A(y)&
    \stackrel{\phi}{\longrightarrow} & p_A(y') \\[2mm]
    \mbox{$\scriptstyle p_A(f)$}
    \bigg
    \uparrow
    \phantom{\mbox{$\scriptstyle p_A(f)$}} & &
    \phantom{\mbox{$\scriptstyle p_A(g)$}}
    \bigg
    \downarrow
    \mbox{$\scriptstyle p_A(g)$} \\ 
    p_A(\widetilde{y})&
     \stackrel{\widetilde{\phi}}{\longrightarrow} 
   & p_A(\widetilde{y}').   \\
   \end{array}\]

By using $p_A(f)_!(\widetilde{y})\simeq y$,
$p_A(g)^*(\widetilde{y}')\simeq y'$,
and fiberwise adjointness,
we have an equivalence
\[ 
    {\rm Map}_{Y_{(a,o)}}(y,
    \phi^*(y'))
    \simeq
    {\rm Map}_{Y_{(\widetilde{a},o)}}(\widetilde{y},
    \widetilde{\phi}^*(\widetilde{y}')),\]
where $a=p_A(y)$ and $\widetilde{a}=p_A(\widetilde{y})$.
Combining equivalences
${\rm Map}_{Y_o}^{\phi}(y,y')\simeq
    {\rm Map}_{Y_{(a,o)}}(y,\phi^*(y'))$
and 
${\rm Map}_{Y_o}^{\widetilde{\phi}}(\widetilde{y},\widetilde{y}') 
    \simeq 
    {\rm Map}_{Y_{(\widetilde{a},o)}}(\widetilde{y},
    \widetilde{\phi}^*(\widetilde{y}'))$,
we obtain the desired equivalence.
\qed

\begin{lemma}\label{lemma:EYA-phi-monoid}
The restriction $B_{\phi}: E(Y/A)_{\phi}\to\mathcal{S}$
is an $E(Y/A)_{\phi}$-monoid object
for any $\phi\in {\rm Tw}^l(A)$.
\end{lemma}

\proof
Suppose that $\phi: s\to t$ in $A$.
If $o\simeq o_1\oplus\cdots\oplus o_n$,
where $o_i\in\mathcal{O}$ for $1\le i\le n$,
then $Y_o\simeq Y_{o_1}\times_A\cdots\times_AY_{o_n}$.
We take $y\in Y_s$ and $y'\in Y_t$
which correspond to
$(y_1,\ldots,y_n)$ and $(y_1',\ldots,y_n')$
under this equivalence, 
respectively.
We have to show that 
$B_{\phi}$ induces an equivalence
${\rm Map}_{Y_o}^{\phi}(y,y')\simeq
   {\rm Map}_{Y_{o_1}}^{\phi}(y_1,y_1')
   \times\cdots\times
   {\rm Map}_{Y_{o_n}}^{\phi}(y_n,y_n')$.
This follows from 
the naturality of fiberwise mapping spaces
and an equivalence 
${\rm Map}_{Y_o}(y,y')\simeq
   {\rm Map}_{Y_{o_1}}(y_1,y_1')
   \times_{{\rm Map}_A(s,t)}\cdots\times_{{\rm Map}_A(s,t)}
   {\rm Map}_{Y_{o_n}}(y_n,y_n')$.
\qed

\bigskip

By Lemmas~\ref{lemma:PO-star-extesion} 
and \ref{lemma:EYA-phi-monoid},
the functor $B$
lies in the essential image
of the fully faithful functor
${\rm Nat}(I,\mathbb{P}_{\mathcal{O}}^*\circ H)_A
\to
  {\rm Map}_{\wcat}
(E(Y/A),\mathcal{S})$
in 
Proposition~\ref{prop:fully-faithful-Nat-Fun},
and hence we obtain a natural transformation
$\beta^A: I\circ {\rm St}(p)\to 
\mathbb{P}_{\mathcal{O}}^*\circ H\circ {\rm St}(p)$.
In particular,
by taking ${\rm St}(p)$ as the identity functor
$A^{\rm op}={\rm Mon}_{\mathcal{O}}^{\rm lax,R}(\cat)\to
{\rm Mon}_{\mathcal{O}}^{\rm lax,R}(\cat)$,
we obtain a natural transformation
$\beta: I\to \mathbb{P}_{\mathcal{O}}^*\circ H$.
We recall that 
$H(\mathcal{C})=\mathcal{C}$
for each $\mathcal{C}\in {\rm Mon}_{\mathcal{O}}^{\rm lax,R}(\cat)$.
By construction,
we see that
$\beta_{\mathcal{C}}\simeq J_{\mathcal{C}}$.
Therefore,
we obtain the following corollary.

\begin{corollary}
The functor $\mathbb{P}_{\mathcal{O}}^*\circ H$
extends to an object
$(\beta: I\to \mathbb{P}_{\mathcal{O}}^*\circ H)$
of $\mathcal{F}$
such that
$\beta_{\mathcal{C}}\simeq J_{\mathcal{C}}$
for each $\mathcal{C}\in
{\rm Mon}_{\mathcal{O}}^{\rm lax,R}(\cat)$.
\end{corollary}

By Lemma~\ref{lemma:characterization-initial-F},
we have a canonical equivalence
between the natural transformations
$(\alpha: I\to \mathbb{P}_{\mathcal{O}!})$
and $(\beta: I\to \mathbb{P}_{\mathcal{O}}^*\circ H)$.
Hence,
we obtain an equivalence
$\gamma: \mathbb{P}_{\mathcal{O}!}
\stackrel{\simeq}{\to}\mathbb{P}_{\mathcal{O}}^*\circ H$
such that $\gamma\circ \alpha\simeq \beta$.
Since $\alpha_{\mathcal{C}}\simeq J_{\mathcal{C}}$
and $\beta_{\mathcal{C}}\simeq J_{\mathcal{C}}$,
we obtain that 
$\gamma_{\mathcal{C}}\simeq {\rm id}_{\mathcal{C}}^*$.
Therefore,
by Lemma~\ref{lemma:characterizetion-eq},
we can regard 
the pair $(H,\gamma)$
as an object of 
${\mathcal Eq}$
and we obtain the following theorem
by Theorem~\ref{thm:uniqueness-equivalences}.

\begin{theorem}\label{thm:T-eq-H}
The pair $(T,\theta)$
is canonically 
equivalent to 
$(H,\gamma)$.
\end{theorem}

This theorem implies the following corollary.

\begin{corollary}\label{cor:main-claim}
The equivalence $T$ in \cite{Torii4}
is equivalent to the restriction $H$
of the equivalence in \cite{HHLN1}.
\end{corollary}




\end{document}